\documentclass[a4paper,12pt]{article}

\usepackage[english]{babel}
\usepackage{graphicx}
\usepackage{amsmath}
\usepackage{float}
\usepackage[top=3.8cm,bottom=3.8cm,left=2.54cm,right=2.54cm]{geometry}

\title{\bf Continuous wavelet transform with the Shannon wavelet from the point of view of hyperbolic partial differential equations\footnotetext{The first author was supported by the grant No.~1391 of the Ministry of Education and Science of the Russian Federation within the basic part of research funding No. 2014/349 assigned to Kursk State University. The second author is grateful for financial support of the BOYSCAST Fellowship 2010-2011 program of the Department of Science and Technology of the Government of India. }}
\author{\sc Eugene B. Postnikov$^1$, Vineet K. Singh$^2$, }

\date{\it $^1$ Department of Theoretical Physics, Kursk State University, Radishcheva st. 33, Kursk 305000, Russia; postnicov@gmail.com\\$^2$ Department of Mathematical Sciences,Indian Institute of Technology, Banaras Hindu University, Varanasi, India; vks1itbhu@gmail.com}

\begin{document}

\maketitle

%%%%% Begin Abstract %%%%%%%%%%%
\begin{abstract}
We identify the result of the continuous wavelet transform with the difference of solutions of two hyperbolic partial differential equations, for which wavelet's shift and scale are considered as independent variables on 2D plane. The characteristic property, which follows from the introduced representation is is the fact
that the transform's values inside the triangle defined by two characteristics
($a=\mathrm{const}$, $b =\mathrm{const}$) and crossecting them slopped line on a scale-shift
plane $(a, b)$ are completely and uniquely defined by the value of transform
and its derivatives along the last mentioned line.
\end{abstract}
%%%%% end %%%%%%%%%%%

%\ams{65T60, 35L52, 35Q94}
%\keywords{Shannon wavelet, PDE, Riemann method.}

%%%% Start %%%%%%
\section{Introduction}

The Shannon wavelets, in both real and complex forms are very important example of wavelet, as from pure mathematical point of view \cite{Novikov2011}, where it is a perfect example for sampling-based multiresolution analysis, as for applications to the signal analysis \cite{Teolis1998,Chui2003}.

Particularly, the standard complex Shannon wavelet defined as
\begin{equation}
\psi{(\xi)}=\mathrm{sinc}{(\xi)}e^{-i2\pi\xi},
\label{eq1}
\end{equation}
where $\mathrm{sinc}(\xi)=\sin(\xi)/\pi\xi$, has a piecewise finite support in Fourier domain, i.e. its spectrum is non-zero only within a finite band.

The corresponding continuous wavelet transform,
\begin{equation}
W(a,b)=C(a) \int_{-\infty}^{+\infty}f(t)\psi^{\ast}(\xi)dt,
\label{def_trans}
\end{equation}
where the asterisk means the complex conjugation, $\xi=(t-b)/a$ defining two output parameters: a positive scaling (dilation) ${a}$ and a location (shift) parameter $b$,  and $C(a)$ is the linear (amplitude) norm factor
determined by the condition:
$$
C(a)=\frac{1}{\int_{-\infty}^{+\infty}|\psi(\xi)|dt},
$$
which is equal to $C(a)=1/\pi a$ for the considered case is applicable for the processing narrowband non-stationary signals. The typical practical examples of these ones belongs to bandpass electrotechnical and electronic systems \cite{Pavlik2005} and, particularly, electrophysilogy \cite{Wang2005,Almanji2011}.

The standard methods of the continuous wavelet transform evaluations deal with its discretization and processing with filter bank implementation, or with the usage Fast Fourier Transform as an intermediate step (due to reducing of the convolution (\ref{def_trans}) to the product in the Fourier space). At the same time, it has been proposed recently the method of calculation based on the representation of the wavelet transform as a Cauchy problem for the proper partial differential equation \cite{Postnikov2006}. It has been done for the example of the standard (reduced) Morlet wavelet and successfully generalized on the case of the exact (admissible) Morlet wavelet being applied to the processing astrophysical data \cite{Postnikov2007}. As well, the solutions of partial differential equations were used for construction of new wavelets \cite{Perel2007} adjusted for an analysis of non-stationary wave propagation in inhomogeneous media.

Thus, the main goal of this work is to prove that the continuous wavelet transform with the complex Shannon wavelet can be considered via solutions of Cauchy problems for partial differential equations and to discuss outlooks, which follows from this fact for this very specific kind of wavelets.

\section{Cauchy Problem for  Partial differential equation based on Shannon Wavelets}
\subsection{Complex Shannon Wavelet and its partial derivatives}

Since the Shannon wavelet is the wavelet function in a strict sense, i.e. satisfies the admissibility condition,
it is required to represent them for the calculation using PDE as a difference of two terms:
$$
\psi{(\xi)}=\frac{\sin(\pi\xi)}{\pi\xi}e^{-i2\pi\xi}=\frac{(e^{i\pi\xi}-e^{-i\pi\xi})}{2i\pi\xi}e^{-i2\pi\xi}=\psi_{1}{(\xi)}-\psi_{2}{(\xi)},
$$
where $\psi_{1}{(\xi)}=ie^{-3\pi i\xi}/2\pi\xi$ and $\psi_{2}{(\xi)}=ie^{-\pi i\xi}/2\pi\xi$.

Such a procedure is analogous to the one, used for evaluation of the wavelet transform with the DoG \cite{Cho1997} and exact Morlet \cite{Postnikov2007} wavelets.

Thus, the transform (\ref{def_trans}) will be represented as a difference of integrals
$$
W(a,b)=\int_{-\infty}^{\infty}f(t)(\psi_{1}^{\ast}(\xi)-\psi_{2}^{\ast}(\xi))\frac{dt}{\pi a}=W_1(a,b)-W_2(a,b),
$$
which have the following explicit forms taken in the Cauchy principal value sense:
\begin{eqnarray}
W_1(a,b)&=&\frac{i}{2\pi}V.P.\int_{-\infty}^{+\infty}f(t)\displaystyle{\frac{e^{3\pi i\frac{t-b}{a}}}{t-b}}dt \label{eq3},\\
W_2(a,b)&=&\frac{i}{2\pi}V.P.\int_{-\infty}^{+\infty}f(t)\displaystyle{\frac{e^{\pi i\frac{t-b}{a}}}{t-b}}dt.\label{eq4}
\end{eqnarray}

No let us consider $a$ and $b$ as independent variables for the desired partial differential equation and search the expressions connecting the corresponding to them partial derivatives. First of all, let us note that the partial derivatives of the expressions (\ref{eq3}) and (\ref{eq4}) provide simple integrals containing $b$ in exponentials only:
\begin{eqnarray}
\frac{\partial W_1}{\partial a}&=&\frac{3}{2 a^2}\int_{-\infty}^{+\infty}f(t)e^{3\pi i\frac{t-b}{a}}dt,\\
\frac{\partial W_2}{\partial a}&=&\frac{1}{2 a^2}\int_{-\infty}^{+\infty}f(t)e^{\pi i\frac{t-b}{a}}dt.
\label{difference}
\end{eqnarray}

Thus, differentiating these expressions with respect to $b$, we get the equalities
\begin{equation}
\frac{\partial^{2}{W_1}}{\partial{a}\partial{b}}+ \frac{3\pi i}{a}\frac{\partial{W_1}}{\partial{a}}=0,
\label{eq7}
\end{equation}
which corresponds to the integral transform (\ref{eq3}) and
\begin{equation}
\frac{\partial^{2}{W_2}}{\partial{a}\partial{b}}+ \frac{\pi i}{a}\frac{\partial{W_2}}{\partial{a}}=0.
\label{eq8}
\end{equation}
corresponding to the integral transform (\ref{eq4}).

\subsection{Solution of Cauchy Problem by the Riemann Method }

Note that both expressions (\ref{eq7}) and (\ref{eq8}) have a form of the following hyperbolic partial differential equation (PDE)
\begin{equation}
\frac{\partial^{2}{u}}{\partial{a}\partial{b}}+ \frac{R}{a}\frac{\partial{u}}{\partial{a}}=0,
\label{eq17}
\end{equation}
where by $u$ and $R$ are $\{W_1, W_2\}$ and $\{3\pi i, \pi i\}$ respectively.

This kind of PDE is well-studied \cite{Koshlyakov}. Particularly, it is known that its characteristics are $a=\mathrm{const}$ $b=\mathrm{const}$, and its unique solution will be found if the initial value is determined along the curve intersecting this two strait lines in one point each one. We take the line $b=-ka+b_{0}$ ($k$, $b_0$ are constant) as this curve. As it will shown below, such a choice provides the sufficient simplification of calculations.

\begin{figure}[H]
\begin{center}
\includegraphics[width=0.6\textwidth]{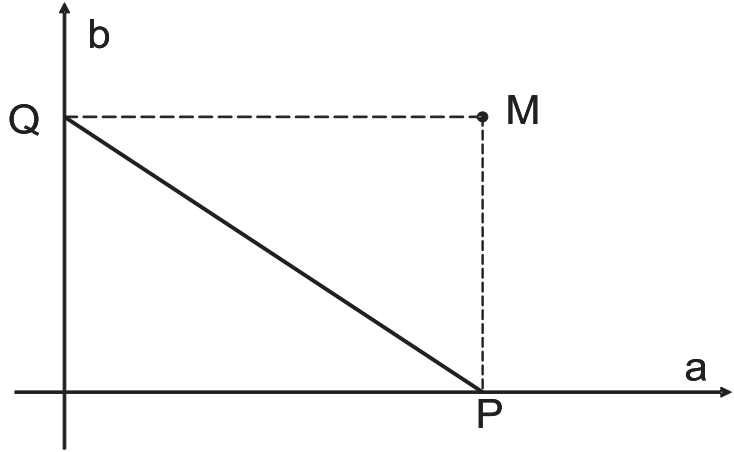}
\end{center}
\caption{The diagram demonstrating the geometry of paths used by Riemann method for the solution of PDE.}
\label{plot}
\end{figure}

The solution approach, suitable for the equation (\ref{eq17}) is called the Riemann method. Its provides the general result \cite{Koshlyakov} for the function $u$ taken in the point $(a_0, b_0)$ (point M in the Fig.~\ref{plot}) in the form
\begin{equation}
u(a_{0},b_{0})=\frac{(uv)_{P}+(uv)_{Q}}{2}+
\frac{1}{2}\int_{QP}\left(v\frac{\partial{u}}{\partial{a}}-u\frac{\partial{v}}{\partial{a}}\right)da-
\left(v\frac{\partial{u}}{\partial{b}}-u\frac{\partial{v}}{\partial{b}}+\frac{2R}{a}u v\right)db,
\label{eq18}
\end{equation}
where points P and Q are denoted in the Fig.~\ref{plot}, and $v(a,b,a_{0}, b_{0})$ is a Riemann function satisfying the  equation conjugated to (\ref{eq17}):
\begin{equation}
\frac{\partial^{2}{v}}{\partial{a}\partial{b}}- \frac{R}{a}\frac{\partial{v}}{\partial{a}}=0
\label{eq19}
\end{equation}

Additionally, according to the general theory, it must also satisfy the following conditions:
\begin{equation}
v(a,b_{0},a_{0},b_{0})=\exp\left[\int_{a_{0}}^{a}0da\right]=1,~~~ \mbox{(on MQ)}
\label{eq20}
\end{equation}
since the Eq. (\ref{eq17}) does not contain a first derivative with respect to $b$, and
\begin{equation}
v(a_{0},b,a_{0},b_{0})=\exp\left[\int_{b_{0}}^{b} \frac{R}{a_{0}}db\right], \quad v(a_{0},b,a_{0},b_{0})=e^{\frac{R(b-b_{0})}{a_{0}}}~~~\mbox{(on MP)}.
\label{eq21}
\end{equation}

It is easy to find the desired general function
\begin{equation}
v(a,b,a_{0},b_{0})=e^{\frac{R(b-b_{0})}{a}},
\label{eq22}
\end{equation}
which satisfies both conditions(\ref{eq20}) and (\ref{eq21}) as well as the equation (\ref{eq19}).

This Riemann function (\ref{eq22}) and its partial derivatives of with respect to $a$ and $b$ along PQ line $b=-ka+b_{0}$ are given as:
\begin{equation}
v(a,b,a_{0},b_{0})=e^{-Rk},
\label{Riemann}
\end{equation}
i.e. it is constant, and
\begin{equation}
\frac{\partial{v}}{\partial{a}}=\frac{-R(b-b_{0})}{a^{2}}e^{\frac{R (b-b_{0})}{a}}=k\frac{R}{a}e^{-Rk},
\label{eq23}
\end{equation}
\begin{equation}
\frac{\partial{v}}{\partial{b}}=\frac{R}{a}e^{\frac{R (b-b_{0})}{a}}=\frac{R}{a}e^{-Rk}.
\label{eq24}
\end{equation}

Since $R=\{3\pi i, \pi i\}$, if we take  $k=2n$, $n\in N^{+}$ then $v=\exp(6n\pi i)=1$ (or $v=\exp(2n\pi i)=1$ as well).
Therefore, the general formula (\ref{eq18}) extremely simplifies:
\begin{equation}
u(a_{0},b_{0})=\frac{(u)_{P}+(u)_{Q}}{2}+
\frac{1}{2}\int_{QP}\left(\frac{\partial{u}}{\partial{a}}-\frac{2Rn}{a}u\right)da-
\left(\frac{\partial{u}}{\partial{b}}+\frac{R}{a}u \right)db.
\label{integral}
\end{equation}

Thus, the final integral representation (\ref{integral}) for the equation (\ref{eq18}) proves that the continuous wavelet transform with the complex Shanon wavelet can be evaluated as a difference of solutions for  Cauchy problems for two hyperbolic partial differential equations (\ref{eq7}), (\ref{eq8}), if the initial values are determined as the transformed function and its partial derivatives along the finite line connecting points on two characteristics of these PDEs.

\section{Example}

Let us consider the simplest harmonic function $f(t)=\exp(i\omega t)$, which has well-known band-limited (with respect to a scale) transform.

\subsection{Initial value for The PDE based on Shannon wavelets}

First of all, we rewrite the equation (\ref{eq3}) in the explicit form:
$$
W_1(a,b)=-\frac{i}{\pi a}V.P.\int_{-\infty}^{+\infty} e^{i\omega t} \frac{e^{3\pi i\xi}}{2\pi\xi}dt=-\frac{i}{2\pi^{2}}e^{-3\pi i \frac{b}{a}}V.P.\int_{-\infty}^{+\infty}\frac{e^{i(\omega +\frac{3\pi}{a})t}}{(t-b)}dt
$$
that gives
\begin{equation}
W_1(a,b)=\frac{1}{2\pi}
\left\{
\begin{array}{rcl}
\exp(i\omega b),& \mathrm{if}& \omega >-3\pi/a,\\
-\exp(i\omega b),& \mathrm{if}& \omega <-3\pi/a.
\end{array}
\right.
\label{eq9}
\end{equation}

One can see that the function (\ref{eq9}) does not depend on $a$, whence the partial derivative $\partial W_1/\partial a \equiv 0$. Another partial derivative is simply
\begin{equation}
\frac{\partial{W_1}}{\partial{b}}=\frac{i\omega}{2\pi}
\left\{
\begin{array}{rcl}
\exp(i\omega b),& \mathrm{if}& \omega >-3\pi/a,\\
-\exp(i\omega b),& \mathrm{if}& \omega <-3\pi/a.
\end{array}
\right.
\label{eq10}
\end{equation}

Analogously, the explicit non-zero results for the transform (\ref{eq4}) are
\begin{equation}
W_2(a,b)=\frac{1}{2\pi}
\left\{
\begin{array}{rcl}
\exp(i\omega b),& \mathrm{if}& \omega >-\pi/a,\\
-\exp(i\omega b),& \mathrm{if}& \omega <-\pi/a.
\end{array}
\right.
\label{eq13}
\end{equation}
and
\begin{equation}
\frac{\partial{W_2}}{\partial{b}}=\frac{i\omega}{2\pi}
\left\{
\begin{array}{rcl}
\exp(i\omega b),& \mathrm{if}& \omega >-\pi/a,\\
-\exp(i\omega b),& \mathrm{if}& \omega <-\pi/a.
\end{array}
\right.
\label{eq14}
\end{equation}

\subsection{Solution of the Cauchy problem}

Let us take the boundary points in such a way that P has co-ordinates $(a_0,0)$ and Q $(a,b_0)$, where $a$ can be arbitrary (within the condition $k=2n$) point admitting $\omega >-3\pi/a$. Then
$$u(P)=\frac{1}{2\pi}~~\mbox{and}~~u(Q)=\frac{e^{i\omega b_{0}}}{2\pi},$$
and substituting (\ref{eq9}), (\ref{eq10})  into (\ref{integral}), we subsequently obtain
$$
W_1(a_{0},b_{0})=\frac{1}{4\pi}+\frac{e^{i\omega b_{0}}}{4\pi}-\frac{1}{2}\int_{QP}\frac{\partial{u}}{\partial{b}}db,
$$
$$
W_1(a_{0},b_{0})=\frac{1}{4\pi}+\frac{e^{i\omega b_{0}}}{4\pi}-\frac{1}{2}\left[u(P)-u(Q)\right],
$$
and, finally,
\begin{equation}
W_1(a_{0},b_{0})=\frac{e^{i\omega b_{0}}}{2\pi},
\label{eq25}
\end{equation}
i.e. the resulting expression (\ref{eq25}) coincides with the initial value (\ref{eq9}) within this subregion, as it must be.

The differential equation (\ref{eq17}) solves by the completely analogous procedure, and its solution coincides with the expression (\ref{eq13}).

Finally, taking the difference between (\ref{eq9}) and (\ref{eq13}), we get the perfect band-like wavelet spectrum, obtained via PDE solutions, the same as a value obtained via the direct transform.

\section{Discussion and outlook}

In this work, we construct the procedure of the continuous wavelet transform with the Shannon wavelet as a solution of a Cauchy problem for the formulated system of partial differential equations. It should by pointed out that this result sufficiently differs from all cases studied early, e.g. \cite{Postnikov2006,Postnikov2007,Perel2007} and other.

This originates from the hyperbolicity of the obtained PDE's. The fact that, due to Riemann's theorem \cite{Koshlyakov}, initial values providing a unique solution, should be stated along a line crosssecting characteristics, which are $a=\mathrm{const}$, $b=\mathrm{const}$ here, do not allow to evaluate this transform using the initial function $f(t)$ taken as the transform at the scale $a=0$.

On the other hand, the transform as a solution of the Cauchy problem exists if the value of transform and its derivatives is known along a finite sloped line on $(a,b)$ plane. Moreover, it follows from the obtained integral representation (\ref{integral}) that this information completely determines all transform's values within the triangle formed by this line and two characteristics. Note also that from the practical point of view, to calculate this one integral with variable limits (P and Q) requires smaller number of operations in comparison with the multiple direct evaluation (for each $a$, $b$ on 2D plane) of the integral (\ref{def_trans}) with the wavelet (\ref{eq1}). As well, we suppose that this result may propose development of evaluation methods for more sophisticated wavelets, e.g. the wave-packet wavelets proposed in \cite{Perel2007}, which constructed on hyperbolic equations too.

\end{document}